\documentclass[a4paper]{amsart}
\pdfoutput=1
\usepackage{amsmath,amssymb,amsfonts,amsthm,longtable}
\usepackage{enumerate}
\usepackage[utf8]{inputenc}
\usepackage{longtable}
\usepackage{a4wide}
\usepackage{graphicx}

\newtheorem{theorem}{Theorem}
\newtheorem{lemma}[theorem]{Lemma}
\newtheorem{proposition}[theorem]{Proposition}

\newtheorem{remark}[theorem]{Remark}
\theoremstyle{definition}

\newtheorem{example}{Example}

\newcommand{\perms}{\mathcal{S}}
\newcommand{\reduct}{\mathsf{red}}
\newcommand{\asc}{\mathsf{asc}}

\newcommand{\invos}{I}

\newcommand{\dgp}{bi-vincular pattern}
\newcommand{\Dgps}{Bi-vincular patterns}
\newcommand{\dgps}{bi-vincular patterns}

\title[Bi-vincular patterns]{Wilf classification of bi-vincular permutation patterns}
\author[R. Parviainen]{Robert Parviainen}
\date{\today}
\thanks{The author acknowledges support by grant no. 090038011 from the Icelandic Research Fund.}
\email{parviainen@ru.is}
\address{The Mathematics Institute, School of Computer Science, Reykjav\'ik University, 103 Reykjav\'ik, Iceland}
\begin{document}
\begin{abstract}
We classify all bi-vincular patterns of length two and three according to the number of permutations avoiding them. These patterns were recently defined by Bousquet-M\'elou et.~al., and are natural generalizations of Babson and Steingr{\'\i}msson's generalized patterns. 
	The patterns are divided into seven and 24 Wilf classes, for lengths two and three, respectively. For most of the patterns an explicit form for the number of permutations avoiding the pattern is given.
\end{abstract}

\maketitle
\section{Introduction}
The study of permutation patterns has, since its relatively slow start, seen an explosive increase in interest with several hundreds
of papers published in the last decade. Even though the applications to other
fields both inside and outside of combinatorics are far reaching, the basic
problem is simple. 

A pattern is a subsequence of the permutation obeying some conditions. The first problem is normally to count the number of permutations that avoid the pattern, that is, that have no subsequence obeying the conditions dictated by the pattern. 

A set of patterns can be \textit{Wilf classified}: divided into equivalence classes where for every pattern in a class the number of permutations of length $n$ avoiding the pattern is the same, for all $n$.  

In this work we Wilf classify all patterns, of lengths up to three, of a new type of pattern. These patterns were recently introduced, \cite{Bousquet-Melou2009}, and have been dubbed \textit{\dgps}.

\subsection{Outline}
After formal definitions are given in the next section, patterns of length two are dealt with in Section \ref{sec:l2}, and of length three in sections \ref{sec:l3a} and \ref{sec:l3b}. Results are summarized in Tables \ref{tab:smalltab} and \ref{tab:bigtab}.


\section{Preliminaries}
Let $\perms_n$ be the set of permutations on $[n]=\{1, 2,\ldots, n\}$. We will use one line notation for permutations, such as $\pi=312$ for the permutation taking 1 to $\pi_1=3$, etc., occasionally adding commas and/or parentheses for clarification, e.g.\ $\pi=(\pi_1,\pi_2,\pi_3)=(3, 1, 2)$. Concatenation of strings will be denoted by $\circ$, so if $a=(1, 2)$ and $b=(4, 3)$ then $a\circ b=(1, 2, 4, 3)$.

For a finite set $S=\{s_1, s_2, \ldots, s_k\}$ with $s_1<s_2<\cdots<s_k$, and a permutation $\sigma$ of $S$, let the \textit{reduction} of $\sigma$, denoted  $\reduct(\sigma)$, be the permutation $\pi$ on $[k]$  such that $\pi_i=j$ if $\sigma_i=s_j$. 

A (classical) pattern is a permutation $p$ on $[k]$. It occurs in $\pi$ if there exists a subsequence $s=(\pi_{i_1}, \pi_{i_2}, \ldots, \pi_{i_k})$ of $\pi$ such that $\reduct(s)=p$. In anticipation of generalized patterns, we write $p$ with dashes, such as $p=$1--3--2, which occurs in $\pi$ if there exist $i<j<k$ such that $\pi_i<\pi_j$, $\pi_i<\pi_k$ and $\pi_j>\pi_k$. 

Generalized patterns, introduced by Babson and Steingr{\'\i}msson, \cite{Babson2000}, add constraints on the elements in the pattern, by specifying elements that must be adjacent. The restrictions are commonly showed by \textit{absence} of dashes. An absence of a dash indicates that the two elements must be adjacent in the permutation as well. For instance, the generalized pattern 13--2 occurs in $\pi$ if there exists $i$ and $j>i+1$ such that $\reduct(\pi_i,\pi_{i+1},\pi_j)=132$.

We can also restrict the first and last elements, by saying that the first (last) element in the pattern must be the first (last) element of the permutation. This may be indicated by bars, so for example $|$3--21$|$ occurs in $\pi\in\perms_n$ if $\pi_1>\pi_{n-1}>\pi_n$.

Bi-vincular patterns were introduced in \cite{Bousquet-Melou2009}, where it was shown that permutations avoiding one particular pattern are in bijection with several other classes of combinatorial objects. The definition uses triples $p=(\sigma, X,Y)$, where $\sigma$ is a permutation on $[k]$ and $X$ and $Y$ are subsets of $\{0\}\cup [k]$.

An occurrence of $p$ in $\pi$ is a subsequence $q=(\pi_{i_1},\ldots,\pi_{i_k})$ such that $\reduct(q)=\sigma$ and, with $(j_1<j_2<\ldots < j_k)$ being the set $\{\pi_{i_1}, \pi_{i_2},\ldots, \pi_{i_k}\}$ ordered (so $j_1=\min_m \pi_{i_m}$ etc.), and $i_0=j_0=0$ and $i_{k+1}=j_{k+1}=n+1$,
\begin{equation*}
	i_{x+1}=i_x+1\ \forall x\in X\mbox{ and } j_{y+1}=j_y+1\ \forall y\in Y.
\end{equation*}
The $X$-set gives constraints on the indices in $q=(\pi_{i_1},\ldots,\pi_{i_k})$, and the $Y$-set on the letters. There are $4^{k+1}k!$ \dgps~ of length $k$, that is, with $\sigma\in\perms_k$.

\begin{example}
A permutation $\pi\in\perms_n$ contains $(132,\{0, 1\},\{2\})$ if it is such that $\pi_1<\pi_2=\pi_{k}+1$ for some $k>2$. An example of such a permutation is $25143$.
\end{example}

\Dgps~of the form $(\sigma,\emptyset,\emptyset)$ are the classical ones, and generalized patterns, also now called \textit{vincular} are those of the form $(\sigma, X,\emptyset)$.

A few instances of these patterns have occurred in the literature, such as irreducible permutations, \cite{Albert2007}, permutations such that $\pi_{i+1}\not=\pi_i-1$. In our notation, they are permutations avoiding $(21,\{1\},\{1\})$. Strongly irreducible permutations, \cite{Atkinson2002}, are those avoiding $(21,\{1\},\{1\})$ and $(12,\{1\},\{1\})$. Permutations avoiding $(231,\{1\},\{1\})$ were studied in \cite{Bousquet-Melou2009}.

For a pattern $p$, let $A_n(p)$ denote the set of permutations in $\perms_n$ that avoids $p$, and let $a_n(p)$ denote $|A_n(p)|$, the number of permutations avoiding it. Also, let $A(p)=(A_1(p), A_2(p), \ldots)$ and $a(p)=(a_1(p), a_2(p), \ldots)$.

Two patterns $p$ and $q$ are said to be \textit{Wilf equivalent} if $a(p)=a(q)$.  

By the \textit{distribution} of a pattern, we mean the two dimensional array $d_{i,j}=$ the number of permutations of length $i$ with exactly $j$ occurrences of the pattern.

We will often use dot diagrams of permutations. These are simply $n$ by $n$ squares of unit boxes, with a dot in the box in column $i$ and row $j$ if and only if $\pi_i=j$. 


\subsection{Symmetries}
Following \cite{Bousquet-Melou2009}, we define inverse, reverse, and complement of \dgps. If $p=(\sigma,X,Y)$, let 
\begin{align*}
	&p^{i}=(\sigma^{i},Y,X),\\
	&p^{r}=(\sigma^{r},k-X,Y),\\
	&p^{c}=(\sigma^{c},X,k-Y),\\
\end{align*}
where $k-S=\{k-s:s\in S\}$ and $\sigma^a$ is the standard \textit{i}nverse, \textit{r}everse, or \textit{c}omplement map of permutations. That is, $\pi^i$ is the permutation such that $\pi^i_k=j$ if $\pi_j=k$, $\pi^r$ is such that $\pi^r_i=\pi_{n+1-i}$ and $\pi^c$ such that $\pi^c_i=n+1-\pi_i$, where $\pi\in\perms_n$. We will also use the more common notation $\pi^{-1}=\pi^i$.

Further, we define compositions of these maps by $p^{ab}=(p^a)^b$. The following results are immediate from the definitions.
\begin{proposition}
	If $p$ occurs in $\pi$, then $p^a$ occurs in $\pi^a$, where $a\in\{i,r,c\}$.
\end{proposition}

Two patterns $p$ and $q$ are in the same \textit{symmetry class} if $q=p^w$ for some finite word $w$ on $\{i,r,c\}$. The 128 \dgps~of length are reduced into 24 symmetry classes, and the 1536 \dgps~ of length 3 into 212 symmetry classes. The number of Wilf (equivalence) classes are seven for length 2, and 24 for length 3.

\begin{example}
The symmetry class of $p=(132,\{0, 1\},\{2\})$ is 
\begin{align*}
\big\{&p=(132,\{0, 1\},\{2\}),
p^i=(132,\{2\},\{0, 1\}),
p^{rc}=(213,\{2, 3\},\{1\}),
p^{irc}=(213,\{1\},\{2, 3\}),\\
&p^r=(231,\{2, 3\},\{2\}),
p^{ir}=(231,\{1\},\{0, 1\}),
p^c=(312,\{0, 1\},\{1\}),
p^{ri}=(312,\{2\},\{2, 3\})\big\}.
\end{align*}
\end{example}
\subsection{The number of symmetry classes}
Burnside's lemma may be used to count the number of symmetry classes. A useful graphical representation of a \dgp~$p=(\sigma, X, Y)$ is the following. Draw the dot diagram for $\sigma\in\perms_n$. Number the horizontal and vertical lines from $0$ to $n$, from bottom to top, and from left to right. Make vertical line $x$ bold if $x\in X$, and horizontal line $y$ bold if $y\in Y$. 
 
With this representation, the symmetries become apparent. Each pattern can be symmetric under one of the following eight transformations: Identity, reflection in the main  diagonal, the anti-diagonal, in a horizontal line, in a vertical line, and rotation by 90, 180, and 270 degrees. To count the the number of patterns fixed by each transformation, we can count separately the number of permutations fixed by them, and the number of configuration of bold lines fixed by them. Let the subscripts $i$, $d$, $a$, $h$, $v$, $90$, $180$, and $270$ indicate the identity transformation, the reflections, and the rotations (with obvious abbreviations). Let $g_x$ and $f_x$ denote the number of permutations and line configurations, respectively, fixed by transformation $x$. 

The number of permutations fixed by the transformations are well known. Of course, $g_i(n)=n!$, and $f_i(n)=4^{n+1}n!$. 

Let $\invos(n)$ be the number of involutions of $[n]$. Then $g_d(n)=g_a(n)=\invos(n)$. 

No permutations are fixed by the other two reflections, so $g_{h}(n)=g_v(n)=0$. 

We also have $g_{270}(n)=g_{90}(n)$, and for $n=4m+k\geq 2$, 
$g_{90}(4m) = g_{90}(4m+1) = \frac{2(2m-1)!} {(m-1)!}$,  and $g_{90}(4m+2) = g_{90}(4m+3) = 0$. 

Finally, $g_{180}(2n) = g_{180}(2n+1) = 2^nn!$.

Counting line configurations fixed by the transformations is easy. In each case, due to the restrictions, we just have to decide which lines should be bold, for a subset of all lines.

Reflection in the main diagonal: Once we chosen which of horizontal lines 1 to $n$ and vertical line $0$ are bold, the remaining configuration is fixed. Thus $f_{d}(n)=f_{a}(n)=2^{n+1}$.

Rotation by 90 degrees: If $n$ is even, we have to pick which of the $n/2+1$ top-most horizontal lines that are bold. If $n$ is odd, which of the top-most $(n+1)/2$ lines. Thus, $f_{90}(n)=f_{270}(n)=2^{\lfloor n/2+1\rfloor}$.

Rotation by 180 degrees: If $n$ is even, we need to choose which of the $n/2+1$ top-most horizontal and $n/2+1$ left-most vertical lines to be bold. If $n$ is odd, which of the top and left $(n+1)/2$ ones.  Thus, $f_{180}(n)=4^{\lfloor n/2+1\rfloor}$.

Now, by Burnside's lemma, the number of symmetry classes $s_n$ of \dgps~ of length $n$ is
\begin{align*}
	s_n&=\frac{1}{8}\left(f_i(n)g_i(n)+f_{180}(n)g_{180}(n)+2f_{90}(n)g_{90}(n)+2f_{d}(n)g_{d}(n)\right).
\end{align*}
Plugging in the expressions for the counts, gives the following result.
\begin{theorem}
The number of symmetry classes $s_n$ of bi-vincular patterns of length $n=4m+k\geq 2$ is given by
\begin{align*}
s_{4m}&=2^{6m-1}(2m)!+2^{8m-1}(4m)!+2^{4m-1}I(4m)+2^{2m}\frac{(2m-1)!}{(m-1)!},\\
s_{4m+1}&=2^{6m-1}(2m)!+2^{8m+1}(4m+1)!+2^{4m}I(4m+1)+2^{2m}\frac{(2m-1)!}{(m-1)!},\\
s_{4m+2}&=2^{6m+2}(2m+1)!+2^{8m+3}(4m+2)!+2^{4m+1}I(4m+2),\\
s_{4m+3}&=2^{6m+2}(2m+1)!+2^{8m+5}(4m+3)!+2^{4m+2}I(4m+3),
\end{align*}
where $\invos(n)$ is the number of involutions on $[n]$,
\[\invos(n)=\sum_{k=0}^{\lfloor n/2\rfloor}\frac{n!}{k!(n - 2 k)!2^k}.\]
\end{theorem}
This confirms that $s_2=24$, $s_3=212$, and gives $s_4=3220$, $s_5=61924$, and $s_7=1478528$.
\subsection{Methods}
Much of the preliminary classification work could be automated, using the software Sage\footnote{\texttt{www.sagemath.org}}. A program was written that finds the symmetry classes, and for a representative from each class enumerates the number of permutations avoiding the pattern up to some small length, looks up the sequence in the Online Encyclopedia of Integer Sequences (OEIS), \cite{Sloane2009}, and writes the results to a file.

\section{Length 2}\label{sec:l2}
We devote a subsection to each Wilf class. The section title gives representatives from each symmetry class, except when they are too large. Also in the title, if it existed at the time of writing, the entry in the OEIS for the sequence of the number of permutations avoiding one of the patterns is given. In the following, A followed by six integers denotes the corresponding sequence in the OEIS. Also, we will let $\mathbf{1}_A$ be equal to 1 if $A$ is true, and 0 if $A$ is false.

The results are summarized in Table \ref{tab:smalltab}.

\begin{table}[ht]
	\centering
	\begin{tabular}{|l|l|l|l|}
		\hline
		Representative & Number of avoiders & Sequence & Section \\
		\hline
		\hline
		$(12,\emptyset,\emptyset)$ & 1, 1, 1, 1, 1, 1, 1 & A000012 &\ref{b01}  \\
		$(12,\emptyset,\{0\})$ & 1, 1, 2, 6, 24, 120, 720 & A000142 &\ref{b02}  \\
		$(12,\{1\},\{1\})$ & 1, 1, 3, 11, 53, 309, 2119 & A000255 &\ref{b03} \\
		$(12,\emptyset,\{0, 1\})$ &1, 1, 3, 12, 60, 360, 2520  & A001710 &\ref{b04} \\
		$(12,\{0\},\{0\})$ & 1, 1, 4, 18, 96, 600, 4320 & A094258 &\ref{b05}  \\
		$(12,\{0, 1\},\{0, 1\})$ & 1, 1, 5, 22, 114, 696, 4920& $n!-(n-2)!$ & \ref{b06} \\
		$(12,X,\{0, 1, 2\})$ & 1, 1, 6, 24, 120, 720, 5040 & $n!-\mathbf{1}_{n=2}$ & \ref{b07} \\
		\hline
	\end{tabular}
	\caption{Wilf classes of \dgps~of length 2. One representative from each class is listed in the first column, and the first seven numbers from the avoidance sequence $a_n(p)$ in the second. If known, a closed form or the sequence entry in the OEIS is given in the third column, and in the last column a reference to the appropriate section is given.}
	\label{tab:smalltab}
\end{table}

\subsection{$(12,\emptyset,\emptyset)$ and $(12,\emptyset,\{1\})$ (A000012)}\label{b01}
All permutations except $(n,n-1,\ldots, 1)$ contain these patterns, so $a_n(p)=1$ for both.


\subsection{$(12,\emptyset,\{0\})$, $(12,\{0\},\{1\})$, and $(12,\{0\},\{2\})$ (A000142)}\label{b02}
It is easy to see that for each pattern, the number of permutations containing at least one occurrence of the pattern is $(n-1)(n-1)!$, since the restrictions are such that one element can be chosen in $(n-1)$ ways, and the remaining elements permuted freely.
\begin{proposition}
	For any pattern $p$ in the class, $a_n(p)=(n-1)!$.
\end{proposition}


\subsection{$(12,\{1\},\{1\})$ (A000255)}\label{b03}
These are the permutations avoiding the substring $(k,k+1)$, which are known to be counted by $a_n=\sum_{i=0}^n (-1)^i (n-i+1) \frac{n!}{i!}$.


\subsection{$(12,\emptyset,\{0, 1\})$ and $(12,\emptyset,\{0, 2\})$ (A001710)}\label{b04}
\begin{proposition}
	For $n>0$ and any pattern $p$ in the Wilf class, $a_n(p)=n!/2$.
\end{proposition}
\begin{proof}
	Assume $p=(12,\emptyset,\{0, 1\})$. We count permutations $\pi$ with at least one occurrence of $p$. Let $i=\pi^{-1}_1$  be fixed. Then there are $(n-i)$ possible choices for $j=\pi^{-1}_2$, and the remaining $n-2$ elements can be chosen without restrictions. Summing over the possible choices of $i$ proves the result.
\end{proof}


\subsection{$(12,\{0\},\{0\})$, $(12,\{0\},\{0, 1\})$, $(12,\{0\},\{0, 2\})$, $(12,\{0\},\{1, 2\})$, $(12,\{1\},\{0, 1\})$, and\newline $(12,\{1\},\{0, 2\})$, (A094258)}\label{b05}
The following proposition is very easily established by counting permutations with at least one occurrence of a pattern.
\begin{proposition}
For any pattern $p$ in the Wilf class, $a_n(p)=n!-(n-1)!$.
\end{proposition}


\subsection{$(12,\{0, 1\},\{0, 1\})$, $(12,\{0, 1\},\{0, 2\})$, $(12,\{0, 1\},\{1, 2\})$, and $(12,\{0, 2\},\{0, 2\})$}\label{b06}
Again, it is easy to enumerate permutations with at least one occurrence of a pattern. 
\begin{proposition}
For $n>1$ and any pattern $p$ in the Wilf class, $a_n(p)=n!-(n-2)!$.
\end{proposition}


\subsection{$(12,X,\{0, 1, 2\})$ and $(21,X,\{0, 1, 2\})$}\label{b07}
It is immediate that the only permutations containing one of the patterns are, respectively, $12$ and $21$.
\begin{proposition}
For any pattern $p$ in the Wilf class, $a_n(p)=n!-\mathbf{1}_{n=2}$.
\end{proposition}


\section{Length 3: Helpful lemmas}\label{sec:l3a}
Before getting to the Wilf classes, we collect a few lemmas that help with further classification.
\begin{lemma}\label{lemma:0Y}
For all pairs $(X,Y)\not\in\big\{\{1\},\{3\},\{0, 1\},\{0, 3\}\big\}\times\big\{\{0, 1\},\{0, 3\}\big\}$, the distributions of $(123,X,\{0\}\cup Y)$ and $(132,X, \{0\}\cup Y)$ are the same.
\end{lemma}
\begin{proof}
For a permutation $\pi$, define a bijection as follows. Reverse the part after the 1, that is, map $\pi=a\circ 1\circ b$ to $\tilde\pi=a\circ 1\circ b^r$. It is easy to check that $\pi$ has an occurrence of $(123,X,\{0\}\cup Y)$ if and only of $\tilde\pi$ has an occurrence of $(132,X,\{0\}\cup Y)$, for all choices of $(X,Y)\not\in\big\{\{1\},\{3\},\{0, 1\},\{0, 3\}\big\}\times\big\{\{0, 1\},\{0, 3\}\big\}$.
\end{proof}

For the next lemma, two slightly different proofs are needed, which is the reason for the dichotomy in the statement.

\begin{lemma}\label{lemma:Y1to3}
Let $A=\big\{\{0\},\{0, 2\}\big\}$ and $B=\big\{\{0, 1, 2\},\{0, 1, 3\},\{0, 2, 3\}\big\}$. For $X\in A\cup B$ and $\sigma\in\{123, 132\}$, the distribution of $(\sigma,X,\{1\})$ is the same as the distribution of $(\sigma,X,\{3\})$.
\end{lemma}
\begin{proof}
For set $A$, we consider the $X=\{0\}$ and $\sigma=123$ case. The exact same proof works for the other cases.

Consider dot diagrams. Occurrences of the patterns are shown in Figure \ref{fig:dotbij1}.
\begin{figure}
\begin{center}
\includegraphics[width=6cm]{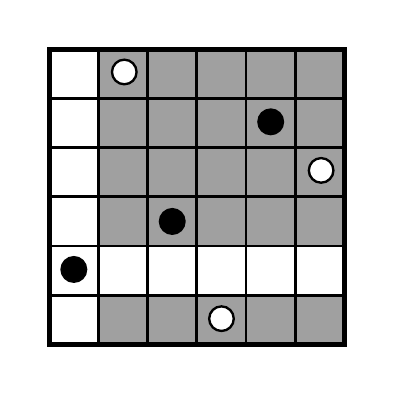}
\quad
\includegraphics[width=6cm]{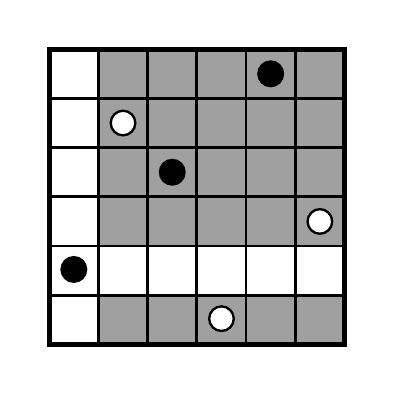}
\end{center}
\caption{On the left, a permutation with an occurrence of $(123,\{0\},\{1\})$, on the right a permutation with an occurrence of $(123,\{0\},\{3\})$. Note that if the shaded rectangles are rotated 180 degrees, the left drawing becomes the right.}
\label{fig:dotbij1}
\end{figure}
A bijection that maps each permutation with $k$ occurrences of $(123,\{0\},\{1\})$ to a permutation with $k$ occurrences of $(123,\{0\},\{3\})$ can be described as follows. In the dot diagram, consider the box with lower left and upper right corners $(1,\pi_1)$ and $(n,n)$ 
(marked with gray background in Figure \ref{fig:dotbij1}). Rotate the box 180 degrees. Simultaneously, rotate the box with upper left and lower right corners $(2,\pi_1{-1})$ and $(n, 0)$ 180 degrees. 

From set $B$, consider the $X=\{0, 1, 2\}$ and $\sigma=132$ case. Again, the exact same proof works for the other cases. Note in particular that in all cases a permutation can have at most one occurrence of the pattern. Example dot diagrams for the $Y=\{1\}$ and $Y=\{3\}$ cases are shown in Figure \ref{fig:dotbij2}.
\begin{figure}
\begin{center}
\includegraphics[width=6cm]{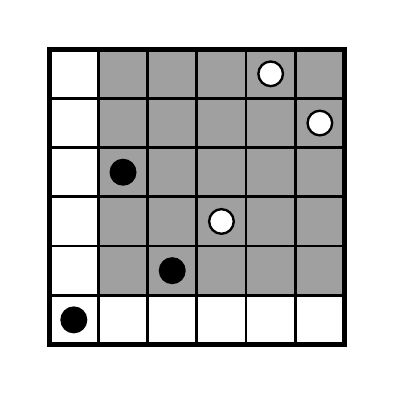}
\quad
\includegraphics[width=6cm]{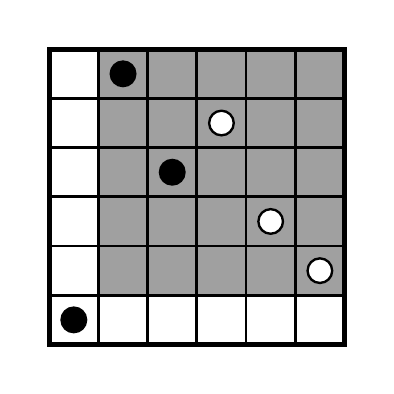}
\end{center}
\caption{In the left drawing, a permutation with an occurrence of $(132,\{0, 1, 2\},\{1\})$, and in the right drawing a permutation with an occurrence of $(132,\{0, 1, 2\},\{3\})$. Notice that if the shaded rectangle is cyclically shifted two steps up, the left drawing becomes the right.}
\label{fig:dotbij2}
\end{figure}
Assume that in the occurrence $(\sigma_1,\sigma_2,\sigma_3)$ of $(132,\{0, 1, 2\},\{1\})$, $\sigma_2=k$.
Cyclically shift all $\pi_i>\pi_1$ $n-k$ steps up, that is, let $\tilde\pi_i=\pi_i+n-k$ if $\pi_1<\pi_i\leq k$, $\tilde\pi_i=\pi_1+\pi_i+k-n$ if $\pi_i>k$, and $\tilde\pi_i=\pi_i$ if $\pi_i<\pi_1$. The shift is illustrated in Figure \ref{fig:dotbij2}.
\end{proof}

We now turn to the Wilf classes for length 3. The results are summarized in Table \ref{tab:bigtab}.

\section{Length 3: Wilf Classes}\label{sec:l3b}
\begin{table}[ht]
	\centering
	\begin{tabular}{|l|l|l|l|}
		\hline
		Representative & Number of avoiders & Sequence & Section \\
		\hline
		\hline
		$(123,\emptyset,\emptyset)$ & 1, 2, 5, 14, 42, 132, 429 & A000108 &\ref{c01}  \\
		$(123,\emptyset,\{1\})$ & 1, 2, 5, 15, 52, 203, 877 & A000110 &\ref{c02}  \\
		$(132,\{1\},\{1\})$ & 1, 2, 5, 15, 53, 217, 1014 & A022493 &\ref{c03} \\
		$(321,\{1\},\{1\})$ & 1, 2, 5, 16, 61, 271, 1372 & A138265 &\ref{c04} \\
		$(132,\{1, 2\},\emptyset)$ & 1, 2, 5, 16, 63, 296, 1623 & A111004 &\ref{c05}  \\
		$(132,\emptyset,\{3\})$ & 1, 2, 5, 16, 64, 312, 1812 & A003149 & \ref{c06} \\
		$(123,\emptyset,\{0\})$ & 1, 2, 5, 16, 65, 326, 1957 & A000522 & \ref{c07} \\
		$(123,\{0\},\{1\})$ & 1, 2, 5, 17, 74, 394, 2484 & A000774 &\ref{c08} \\
		$(123,\{0, 2\},\emptyset)$ & 1, 2, 5, 17, 75, 407, 2619 &  See main text. & \ref{c09} \\
		$(123,\{2\},\{2\})$ & 1, 2, 5, 18, 82, 459, 3041 &  & \ref{c10} \\
		$(123,\emptyset,\{0, 1\})$ & 1, 2, 5, 18, 84, 480, 3240 & $n!-(n-1)!(n-2)/2$ & \ref{c11} \\
		$(132,\{2\},\{1, 2\})$ & 1, 2, 5, 18, 85, 494, 3389 &  & \ref{c22} \\
		$(132,\{1\},\{1, 2\})$ & 1, 2, 5, 18, 86, 502, 3444 &  & \ref{c23} \\
		$(123,\{1\},\{1, 2\})$ & 1, 2, 5, 19, 90, 523, 3573 &  & \ref{c24} \\
		$(123,\{1\},\{1, 3\})$ & 1, 2, 5, 19, 91, 531, 3641 & A052169 & \ref{c12} \\
		$(123,\{0\},\{0\})$ & 1, 2, 5, 19, 97, 601, 4321 & $n!-(n - 1)! + 1$ & \ref{c13} \\
		$(123,\emptyset,\{0, 1, 2\})$ & 1, 2, 5, 20, 100, 600, 4200 & $5n!/6$ & \ref{c14} \\
		$(123,\{0\},\{0, 1\})$ & 1, 2, 5, 20, 102, 624, 4440 & $n! - (n - 2)! (n - 2)$ & \ref{c15} \\
		$(132,\{1, 2\},\{1, 2\})$ & 1, 2, 5, 20, 102, 626, 4458 & See main text. & \ref{c16} \\
		$(123,\{1, 2\},\{1, 2\})$ & 1, 2, 5, 21, 106, 643, 4547 & A002628 & \ref{c17}\\
		$(123,\{0\},\{0, 1, 2\})$ & 1, 2, 5, 21, 108, 660, 4680 & $n! - (n - 1)!/2$ & \ref{c18} \\
		$(123,\{0, 1\},\{0, 1\})$ & 1, 2, 5, 22, 114, 696, 4920 & $n!-(n-2)!$ & \ref{c19} \\
		$(123,\{0, 1, 2\},\{0, 1, 2\})$ & 1, 2, 5, 23, 118, 714, 5016 & $n!-(n-3)!$ & \ref{c20} \\
		$(123,X,\{0, 1, 2, 3\})$ & 1, 2, 5, 24, 120, 720, 5040 & $n!-\mathbf{1}_{n=3}$ & \ref{c21} \\
		\hline
	\end{tabular}
	\caption{Wilf classes of \dgps~of length 3. One representative from each class is listed in the first column, and the first seven numbers from the avoidance sequence $a_n(p)$ in the second. If known, a closed form or the sequence entry in the OEIS is given in the third column, and in the last column a reference to the appropriate section is given.}
	\label{tab:bigtab}
\end{table}

\subsection{$(123,\emptyset,\emptyset)$, $(132,\emptyset,\emptyset)$, and $(132,\emptyset,\{1\})$ (A000108)}\label{c01} 

The first two are classical cases, and part of the Catalan family. Claesson showed in \cite{Claesson2001} that permutations avoiding the third pattern are also counted by the Catalan numbers.


\subsection{$(123,\emptyset,\{1\})$ and $(132,\emptyset,\{2\})$ (A000110)} \label{c02}

These cases were also treated by Claesson, \cite{Claesson2001}, and the number of permutations avoiding one of the patterns are counted by the Bell numbers. 


\subsection*{Interlude on ascent sequences}\label{sec:ascent}

For the next two Wilf classes, some background material from \cite{Bousquet-Melou2009} will be needed.

An \textit{ascent sequence} is a sequence $x_i$, $1\leq i\leq n$ such that $x_1=0$, and for $i>1$, $0\leq x_i\leq 1+\asc_{i-1}$, where $\asc_i$ are the number of ascents (positions $j$ such that $x_j<x_{j+1}$) in $(x_1,\ldots, x_i)$.

In \cite{Bousquet-Melou2009} a bijection $f$ from ascent sequences to $(231,\{1\},\{1\})$-avoiding permutations is given. To define it, we first define \textit{modified} ascent sequences.

Let $x$ be an ascent sequence. Define $\mathcal{A}(x)=\{i: i \in [n-1] \mbox{ and } x_i< x_{i+1}\}$.
Denote by $\hat{x}$ the outcome of the following algorithm.
\medskip\\
\noindent
\begin{minipage}{30em}
\mbox{}{\tt for} $i\in \mathcal{A}(x)$: \\
\mbox{}\qquad{\tt for} $j \in [i-1]$: \\
\mbox{}\qquad\qquad{\tt if} $x_j \geq x_{i+1}$ {\tt then} $x_j := x_j+1$
\end{minipage}
\medskip\\
and call $\hat{x}$ the {\it{modified ascent sequence}}.
For example, if $x=(0,1,0,1,3,1,1, 2)$ then $\mathcal{A}(x)=(1,3,4,7)$ and $\hat{x}=(0,3,0,1,4,1,1, 2)$.

The map $f$ is defined as follows. Let $x$ be an ascent sequence. Write  the two line word with columns $(x_i,i)^\top$. Now sort the columns with respect to top entry, and break ties by sorting in descending order with respect to the bottom entry. The bottom row is $f(x)$. 

The inverse is pretty easy to describe, once we have defined \textit{active sites}. 
For a permutation $\pi\in A_n(p)$, the active sites are the positions $n+1$ can be inserted into to give a permutation in $A_{n+1}(p)$. If $n+1$ is inserted in position $i$, then $\pi\to (\pi_1,\ldots, \pi_i,n+1,\pi_{i+1},\ldots,\pi_n)$. 

For a permutation $\pi$, let $\pi^{(k)}$ be the reduction of the permutation obtained by removing all $\pi_i>k$.

Encode a permutation in $A_n(231,\{1\},\{1\})$ by its insertion word $x=(x_1,x_2, \ldots, x_n)$, where $x_i=j$ if $\pi^{(i)}$ is obtained from $\pi^{(i-1)}$ by inserting $i$ the $j$th active site. It is shown in \cite{Bousquet-Melou2009} that the insertion word equals $f(\pi)$.


\subsection{$(132,\{1\},\{1\})$ and $(231,\{1\},\{1\})$ (A022493)} \label{c03}

The pattern $(231,\{1\},\{1\})$ is studied in \cite{Bousquet-Melou2009}. A central result is that permutations avoiding the pattern are in bijection with $(2+2)$-free posets. The sequence A022493 also counts upper triangular non-negative integer matrices with sum $n$ and no zero rows or columns, see \cite{Dukes2009}.

Wilf equivalence of the two patterns is the next result.

\begin{lemma}
The two patterns $(132,\{1\},\{1\})$ and $(231,\{1\},\{1\})$ are Wilf equivalent.
\end{lemma}
\begin{proof}
Only minor modifications are necessary to turn $f$ into a bijection from $A_n(132,\{1\},\{1\})$ to ascent sequences of length $n$. Combining these two bijections proves the lemma.

Namely, the insertion word of $\pi\in A_n(132,\{1\},\{1\})$ also an ascent sequence. The proof of this is completely analogous to the proof of Theorem 1, \cite{Bousquet-Melou2009}, and thus the map $g(\pi)=x$, where $x$ is the insertion word, is a bijection from $A_n(132,\{1\},\{1\})$ to the set of ascent sequences of length $n$. Therefore, $f(g(\pi))$ is a bijection from $A_n(132,\{1\},\{1\})$  to $A_n(231,\{1\},\{1\})$.
\end{proof}

\begin{remark}
The inverse of $g$ has an easy description. Consider the same two line word as in the definition of $f$. Now, divide the two line word by the non-decreasing runs in the top row. Reverse the word, keeping these runs intact. The bottom row of the result is the permutation.
\end{remark}


\subsection{$(321,\{1\},\{1\})$ (A138265)}\label{c04}
The sequence counts the number of upper triangular binary matrices with $n$ 1's with no zero rows or columns. 

\begin{proposition}
The number of permutations of $[n]$ avoiding $(321,\{1\},\{1\})$ equals the number of upper triangular binary matrices with $n+1$ 1's with no zero rows or columns, and also equals the number of permutations in $\perms_{n+1}$ avoiding both $(231,\{1\},\{1\})$ and $(21,\{1\},\{1\})$.
\end{proposition}
\begin{remark} Compare with the previous class, where the same matrices without the binary restriction are in bijection with $(231,\{1\},\{1\})$-avoiding permutations. 
\end{remark}
\begin{proof}
It is shown in \cite{Dukes2009} that the binary matrices can be represented as ascent sequences, with no two adjacent elements equal. 

The map $f$, restricted to ascent sequences with no two adjacent elements equal, becomes a bijection to permutations avoiding both $(231,\{1\},\{1\})$ \textit{and} $(21,\{1\},\{1\})$.

We define another map $h$, a simple modification of $f$, from these ascent sequences to the wanted permutations. 

For $i=2,\ldots,n$, write the two line word with columns $(x_i,i-1)^\top$. Sort the columns with respect to top entry, and break ties by sorting in \textit{ascending} order with respect to the bottom entry. The bottom row is $h(x)$. Mimicking the proof of Theorem 1, \cite{Bousquet-Melou2009}, it follows that $h$ is a bijection, from ascent sequences of length $n$ with no two consecutive elements equal, to $A_{n-1}(321,\{1\},\{1\})$.
\end{proof}


\subsection{$(132,\{1, 2\},\emptyset)$ (A111004)} \label{c05}

This is the consecutive pattern 132, studied in \cite{Elizalde2003}.


\subsection{$(132,\emptyset,\{3\})$ and $(132,\emptyset,\{1, 3\})$ (A003149)}\label{c06}

\begin{proposition}
$a_n(132,\emptyset,\{3\})=a_n(132,\emptyset,\{1, 3\})=\sum_{k=1}^{n} (k-1)!(n-k)!$.
\end{proposition}
\begin{proof}
Let $k=\pi^{-1}_n$. To avoid either pattern, $\pi_1,\pi_2, \ldots, \pi_{k-1}$ must all be greater than $k$, and $\pi_{k+1},\ldots, \pi_{n}$ all less than $k$. The number of such permutations is $(k-1)!(n-k)$!. Summing over the possible values of $k=\pi^{-1}_n$ gives the result.
\end{proof}


\subsection{$(123,\emptyset,\{0\})$, $(123,\{0\},\{2\})$, $(132,\emptyset,\{0\})$, and $(132,\{0\},\{2\})$ (A000522)} \label{c07}

Equality of the first two with the latter two follows from Lemma \ref{lemma:0Y}. 
Using symmetry, equality of the first two and of the last two follows from the next lemma.

\begin{lemma}
For $\sigma\in\{123, 132\}$, a permutation avoids $(\sigma,\{2\},\{0\})$ if and only if it avoids $(\sigma,\emptyset,\{0\})$.
\end{lemma}

\begin{proof}
We show the case $\sigma=123$. The other case is similar. If a permutation has an occurrence of $(123,\{2\},\{0\})$, it obviously has an occurrence of $(123,\emptyset,\{0\})$.

Assume a permutation has at least one occurrence of $(123,\emptyset,\{0\})$ and consider the part after the 1. Let $k$ be the smallest index such that $\pi_k$ is the element playing the role of ``3'' in an occurrence of $(123,\emptyset,\{0\})$. Then $1<\pi_{k-1}<\pi_k$, because otherwise $\pi_{k-1}$ would be the ``3'' in an occurrence of $(123,\emptyset,\{0\})$.
\end{proof}

\begin{proposition}
The number of permutations of $[n]$ avoiding one of the patterns in this class is given by $a_n=\sum_{k=0}^{n-1}\frac{(n-1)!}{k!}$.
\end{proposition}

\begin{proof}
Consider the pattern $(123,\emptyset,\{0\})$. Assume $i=\pi^{-1}_1$. All letters $\pi_j$, $j>i$, must be in decreasing order, with no restrictions on the other elements. The number of such arrangements are $(n-1)!/(n-i)!$. Sum over the possible values of $\pi_1$.
\end{proof}


\subsection{$(123,\{0\},\{1\})$,
$(123,\{0\},\{3\})$,
$(132,\{0\},\{1\})$,
$(132,\{0\},\{3\})$,
$(132,\{1\},\{3\})$,
$(132,\{2\},\{3\})$, and 
$(132,\{3\},\{3\})$ (A000744).}\label{c08}
\begin{figure}
\begin{center}
\includegraphics[width=14cm]{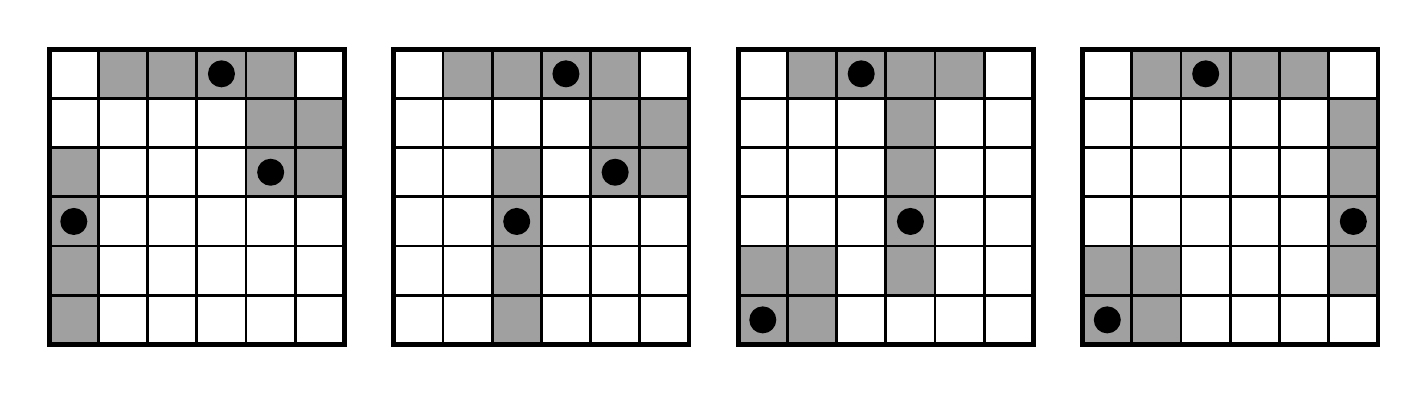}
\caption{Permutations with occurrences of the pattern, from left to right, $(132,\{0\},\{3\})$,
$(132,\{1\},\{3\})$,
$(132,\{2\},\{3\})$, and
$(132,\{3\},\{3\})$.}
\label{fig:bij132}
\end{center}
\end{figure}
Equality of the first with the next three follows from lemmas \ref{lemma:0Y} and \ref{lemma:Y1to3}. Equality of the latter three with the fourth can be deduced from the dot diagrams: in Figure \ref{fig:bij132}, occurrences of the patterns are shown. In each case, the number of permutations avoiding the pattern can be enumerated in the same way as in the proof below.
\begin{proposition}
The number of permutations of $[n]$ that avoid any pattern $p$ in the class is 
\[a_n(p)=(n-1)!(1+\sum_{k=1}^{n-1}\frac{1}{k}).\]
\end{proposition}
\begin{proof}
Consider the pattern $(132,\{0\},\{3\})$. We will enumerate permutations avoiding it. First, choose $\pi_1$ and then $j=\pi^{-1}_n$, and finally consider the possible permutations of the remaining letters. The number of possibilities are
\begin{align*}
&(n-2)!&(\pi_1=1)\\
+&(n-2)!+(n-3)!&(\pi_1=2)\\
+&(n-2)!+2(n-3)!+1\cdot 2(n-4)!&(\pi_1=3)\\
+&(n-2)!+3(n-3)!+2\cdot 3(n-4)!+1\cdot 2\cdot 3(n-5)!&(\pi_1=4)\\
+&(n-2)!+4(n-3)!+3\cdot 4(n-4)!+2\cdot 3\cdot 4(n-5)!+1\cdot 2\cdot 3\cdot 4(n-6)!&(\pi_1=5)\\
+&\cdots\\
+&(n-2)!+(n-2)(n-3)!+(n-3)(n-2)(n-4)!+\cdots+(n-2)!&(\pi_1=n-2)
\end{align*}
Some simple sums finish the proof.
\end{proof}


\subsection{$(123,\{0, 2\},\emptyset)$ and $(132,\{0, 2\},\emptyset)$} \label{c09}

Lemma \ref{lemma:0Y} gives equivalence of the two patterns. 

\begin{proposition}\label{prop:conj}
Let $a_{n,k}=|\{\pi\in A_n(123,\{0, 2\},\emptyset): \pi_1=k\}|$. Then
\[a_{n,k}=\begin{cases}(k-1)!k^{n-k}\quad&\text{for }k=1,2,\ldots, n-2 \\ 
(n-1)!&\text{if }k=n-1\text{ or }k=n 
\end{cases}.\]
\end{proposition}
\begin{proof}
Assume a permutation $\pi\in S_n$ begins with $k$, and that it avoids the pattern $(123,\{0, 2\},\emptyset)$. First note that if $k=n-1$, the number of such permutations are $(n-1)!=(k-1)!k$. Assume now that $k\leq n-2$, and consider the active sites in $\pi$. They are $\pi^{-1}_1,\pi_2^{-1},\ldots, \pi_k^{-1}$, for if and only if $n+1$ is inserted after an element greater than $k$ an occurrence of $(123,\{0, 2\},\emptyset)$ is created. Further, if $n+1$ is removed from $A_{n+1}(123,\{0, 2\},\emptyset)$, the result is a permutation in $A_n(123,\{0, 2\},\emptyset)$. Hence, for $k\leq n-2$, $a_{n+1,k}=ka_{n,k}=(k-1)!k^{n-k}$. 
\end{proof}
At the time of writing, the sequence A018927, which is $(n+1)!-a_{n+1}(123,\emptyset,\{0, 2\})$ in our notation, has as its main description the sum over all permutations of $[n]$ of the maximum of $\pi_i-i$. The formula given is a sum which differs from that achieved by summing over $k$ in Proposition \ref{prop:conj}. Equality of the two sums is the next result.  
\begin{proposition}
We also have
\[a_n(123,\{0, 2\},\emptyset)=n!-\sum_{k=1}^{n-2}k\cdot k!\big((k+1)^{n-k-1}-k^{n-k-1}\big).\]
\end{proposition}

A bijective proof of the equality would be nice, but has eluded the author. An elementary proof starting with the expressions in Proposition  \ref{prop:conj} is not hard to give, but is neither short nor particularly interesting, and hence omitted. For a proof of the interpretation of the sum as the sum of the maximal excedence over all permutations, see \cite[p.~107]{Knuth1998}.


\subsection{$(123,\{2\},\{2\})$ and $(132,\{2\},\{2\})$ (1, 2, 5, 18, 82, 459, 3041,\ldots).} \label{c10}
\begin{figure}
\begin{center}
\includegraphics[width=12cm]{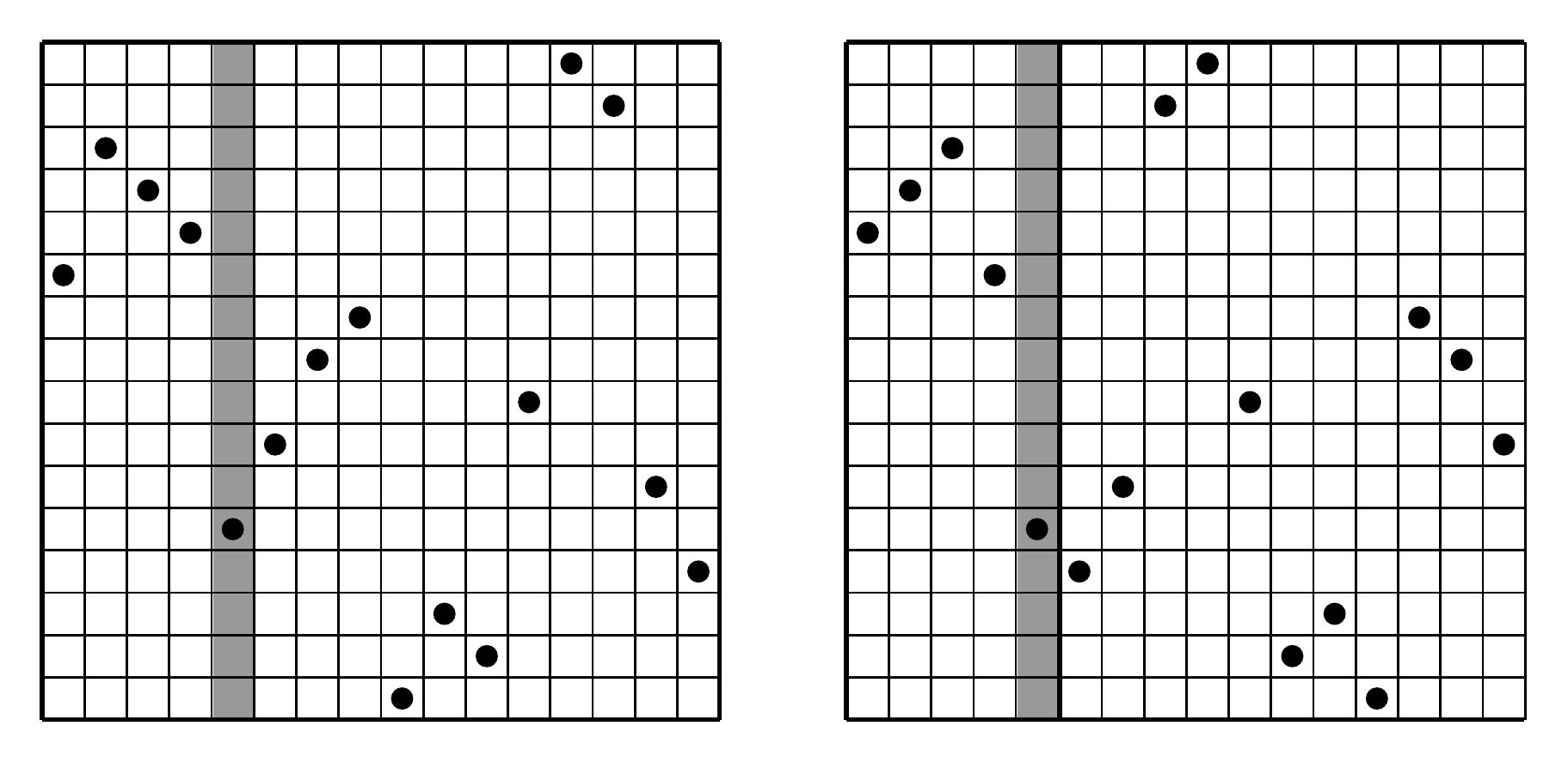}
\end{center}
\caption{Illustration of the bijection used in the proof of Lemma\ref{lemma:wilf22}. On the left, a permutation with exactly one occurrence of $(123,\{2\},\{2\})$-avoiding permutation, and on the right, with one occurrence of $(132,\{2\},\{2\})$.}\label{fig:reflectex}
\end{figure}

\begin{lemma}\label{lemma:wilf22}
The patterns $(123,\{2\},\{2\})$ and $(132,\{2\},\{2\})$ are Wilf equivalent.
\end{lemma}
\begin{proof}
We show that the same recursion holds for $a_n(123,\{2\},\{2\})$ and $a_n(132,\{2\},\{2\})$, namely $a_2=2, a_3=5$, and for $n\geq 4$, 
\[a_{n+1}=c_n+n(a_n-a_{n-1})+(n+1)a_{n-1},\]
where $c_n$ is to be determined.

Let $p_1=(123,\{2\},\{2\})$ and $p_2=(132,\{2\},\{2\})$.

Consider the active sites, as defined in Section \ref{sec:ascent}, of a permutation $\pi\in\perms_n$ avoiding $p_1$. If $\pi$ starts with $n$, then all $n+1$ sites are active, and there are $a_{n-1}(p_1)$ such permutations.

If $\pi$ does not start with $n$, there is exactly one non-active site, after $\pi_i=n$, and there are $a_n(p_1)-a_{n-1}(p_1)$ such permutations. 

However, not all of $A_{n+1}(p_1)$ can be obtained by insertions in the permutations in $A_n(p_1)$. If a permutation has exactly one occurrence of $(12,\{1\},\{1\})$, say $(\pi_k,\pi_{k+1}=\pi_k+1)$, with at least one $\pi_i<\pi_k$, for $i<k$, then the permutation $(\pi_1,\ldots, \pi_k,n,\pi_{k+1},\ldots, \pi_n)$ is in $A_{n+1}(p_1)$. Let $c_n(p_1)$ denote the number of such permutations. Hence,
\[a_{n+1}(p_1)=c_n(p_1)+n\left(a_n(p_1)-a_{n-1}(p_1)\right)+(n+1)a_{n-1}(p_1).\]

Now consider permutations avoiding $p_2$. By the same reasoning as above, it can be seen that 
\[a_{n+1}(p_2)=c_n(p_2)+n\left(a_n(p_2)-a_{n-1}(p_2)\right)+(n+1)a_{n-1}(p_2),\]
where $c_n(p_2)$ is the number of permutations $\pi\in\perms_n$ with exactly one occurrence of $(21,\{1\},\{1\})$, say $(\pi_k,\pi_{k+1}=\pi_k-1)$, such that there exists an $i<k$ with $\pi_i<\pi_k$.

As it is trivial to check the initial values, it only remains to show that $c_n(p_1)=c_n(p_2)$. A suitable bijection with confirm this claim. 

Assume that $\pi$ is a permutation counted by $c_n(p_1)$, that is, there exists a unique $k$ such that  $\pi_k=\pi_{k+1}-1$, and for some $i<k$, $\pi_i<\pi_k$. For $j<k$, let $\pi_j$ be the smallest value that is less than $\pi_k$. Write $\pi=\pi_a\circ\pi_j\circ\pi_b$, and define $h(\pi)=\pi^r_a\circ\pi_j\circ\pi^r_b$. 

The operation is illustrated in Figure  \ref{fig:reflectex}. By construction, the resulting permutation $h(\pi)$ is one counted by $c_n(p_2)$.

The inverse map is easily described, and $h(\pi)$ is  the required bijection.
\end{proof}

\subsection{$(123,\emptyset,\{0, 1\})$ and many more}\label{c11}
See Section \ref{sec:classA} for a list of the class.
\begin{proposition}
Let $b_n$ be the number of permutations of $[n]$ that has at least one occurrence of $(123,\emptyset,\{0, 1\})$. Then $b_n=\frac{1}{2}(n-1)!(n-2)$, the Lah numbers.
\end{proposition}
\begin{remark}
The sequence $b_n$ above is A001286.
\end{remark}
\begin{proof}
Count permutations with at least one occurrence of$(123,\emptyset,\{0, 1\})$. First choose $i=\pi^{-1}_1$ from $\{1,\ldots,n-2\}$, then choose $j=\pi^{-1}_2$ from $\{i+1,\ldots,n-1\}$. The remaining $n-2$ elements can be permuted at will. This gives
\[b_n=(n-2)!\left((n-2)+(n-3)+\cdots +1\right)=\frac{1}{2}(n-1)!(n-2).\]
\end{proof}
Permutations with at least one occurrence of any given pattern in the class can be enumerated in basically the same way as above.
\begin{lemma}
	All patterns in the class are Wilf equivalent.
\end{lemma}

\subsection{$(132,\{2\},\{1, 2\})$ (1, 2, 5, 18, 85, 494, 3389, \ldots).} \label{c22}
No significant results to report.

\subsection{$(132,\{1\},\{1, 2\})$ (1, 2, 5, 18, 86, 502, 3444, \ldots).}\label{c23}
No significant results to report.


\subsection{$(123,\{1\},\{1, 2\})$ (1, 2, 5, 19, 90, 523, 3573, \ldots).} \label{c24}
No significant results to report.


\subsection{$(123,\{1\},\{1, 3\})$ and $(132,\{2\},\{0, 2\})$ (A052169)}\label{c12}

Equality follows from  Lemma \ref{lemma:0Y} (in the symmetry class of the first pattern is $(123,\{2\},\{0, 2\})$).

\begin{proposition}
Let $a_n=a_n(123,\{1\},\{1, 3\})$. Then $a_1=1, a_2=2$, and for $n>2$, \[a_n=(n-1)a_{n-1}+(n-2)a_{n-2}.\] 
The recursion is solved by
\begin{align*}
	a_n&=(n-1)!+\sum_{k=0}^{n-2}(-1)^{n-k}(k+1)!\binom{n-1}{k}\\
		&=\frac{c_{n+1}}{n}=c_n+c_{n-1},
\end{align*}
where $c_n$ is the number of non-derangements of $[n]$,
\[
	c_n=n!(1-\sum_{k=1}^n(-1)^k\frac{1}{k!}).
\]
\end{proposition}
The result follows from the next lemma (after summing over $k$).
\begin{lemma}
Let $a_{n,k}$ be the number of permutations of $[n]$ avoiding $(123,\{1\},\{1, 3\})$ for which $k=\pi_n$. Then $a_{n, 1}=a_{n, 2}=(n-1)!$ and for $k=3,\ldots, n$, $a_{n,k}=a_{n,k-1}-a_{n-1,k-1}$.
\end{lemma}
\begin{proof}
Let $S_{n,k}\subset\perms_n$ be the permutations $\pi$ that avoid $(123,\{1\},\{1, 3\})$, and for which $k=\pi^{-1}_n$. We will give a bijection from $S_{n,k}$ to $S_{n,k-1}\setminus S_{n-1,k-1}$.

Consider again the dot diagrams of the permutations. For each permutation $\pi$ in $S_{n,k}$, swap columns $k-1$ and $k$ in the diagram, and let $\tilde{\pi}$ be the resulting permutation. Let $S'_{n,k}$ be the set of all permutations $\tilde\pi$ such that $\tilde\pi=\pi$ for some $\pi\in S_{n,k}$, and let $T_{n,k-1}=S_{n,k-1}\setminus S'_{n,k}$.

Claim 1: $\tilde\pi\in S_{n,k-1}$.

For all permutations $\pi$ in $T_{n,k-1}$, let $\bar\pi$ be the permutation that results if row $\pi_{k+1}$ and column $k+1$ is removed. Let $T'_{n,k-1}$ be the set of permutations so achieved.

Claim 2: $T'_{n,k-1}=S_{n-1,k-1}$.

Once claims 1 and 2 are established, it follows that the map $f$  defined by $f(\pi)=\tilde\pi$ is the needed bijection. 

To see that claim 1 holds, we just need to check that $\tilde\pi$ avoids $(123,\{1\},\{1, 3\})$, but this is evident since the column swap can not create an occurrence of the pattern (recall that the ``3'' in the pattern must be $n$).

Similarly for claim 2, it is clear that $\bar\pi\in S_{n-1,k-1}$. It still remains to show that for $\sigma\in S_{n-1,k-1}$, there exists a $\pi\in T_{n,k-1}$ such that $\bar\pi=\sigma$.

Consider the permutations in $T_{n,k-1}$. These are the permutations such that if columns $k$ and $k+1$ in the dot diagram are swapped, they have an occurrence of $(123,\{1\},\{1, 3\})$ that consists of $(\pi_{k-2},\pi_{k-1},\pi_{k})$. For a permutation in $S_{n-1,k-1}$, insert a column and row with a dot between rows $k-2$ and $k-1$ and columns $\pi_{k-2}$ and $\pi_{k-2}+1$. The result is a permutation in $T_{n,k-1}$.
\end{proof}


\subsection{$(123,\{0\},\{0\})$, $(123,\{0\},\{0, 2\})$, $(132,\{0\},\{0\})$, and  $(132,\{0\},\{0, 2\})$}\label{c13}
Equivalence of the first and third, and of the second and fourth, follows from Lemma \ref{lemma:0Y}. 

\begin{lemma}
A permutation avoids $(123,\{0\},\{0, 2\})$ if and only if it avoids $(123,\{0\},\{0\})$. 
\end{lemma}
\begin{proof}
We need to show the if part, so assume a permutation has an occurrence, say $(\pi_1,\pi_i,\pi_j)$, of  $(123,\{0\},\{0\})$. If $\pi_i=\pi_j-1$ we are done. Otherwise, consider the elements $\pi_{m_k}$, $1\leq k\leq \pi_j-\pi_i-1=a$, with $\pi_i<\pi_{m_k}<\pi_j$. 

Either $m_k<i$ for all $k$, in which case $(\pi_1,\pi_{a},\pi_j)$ is an occurrence of $(123,\{0\},\{0, 2\})$, or at least one $m_k=b>i$, in which case $(\pi_1,\pi_{b-1},\pi_b)$ is an occurrence of $(123,\{0\},\{0, 2\})$.
\end{proof}
\begin{proposition}
For any pattern $p$ in the class,
\[a_n(p)=n! - (n - 1)! + 1.\]
\end{proposition}
\begin{proof}
If $\pi$ has an occurrence of $p=(123,\{0\},\{0\})$, then $\pi_1=1$. The permutations with at least one occurrence of $p$ are thus all permutations that start with 1, except $(1,n,(n-1),\ldots,2)$. 
\end{proof}
\begin{remark}
The sequence $c_n=(n-1)!-1$ is A033312.
\end{remark}


\subsection{$(123,\emptyset,\{0, 1, 2\})$,
$(123,\emptyset,\{0, 1, 3\})$,
$(132,\emptyset,\{0, 1, 2\})$,
$(132,\emptyset,\{0, 1, 3\})$,
$(132,\emptyset,\{0, 2, 3\})$, and 
$(132,\emptyset,\{1, 2, 3\})$}\label{c14}
\begin{proposition}
For any pattern $p$ in the class,
\[a_n(p)=\frac{5}{6}n!.\]
\end{proposition}
\begin{proof}
We count permutations with at least one occurrence of $(123,\emptyset,\{0, 1, 2\})$. There exist $i<j<k$ such that  $\pi_i=1$, $\pi_j=2$ and $\pi_k=3$. Choose $i$, then $j$, then $k$, then the remaining $n-3$ elements freely. The number of possibilities are
\begin{align*}
\Big(&(n-1)\big((n-2)+(n-3)+\cdots+1\big)+\\&(n-2)\big((n-3)+(n-4)+\cdots+1\big)+\cdots+ 1\Big)(n-3)!=\frac{n!}{6}.
\end{align*}
The other six cases are very similar.
\end{proof}
\begin{remark}
The sequence $n!-5n!/6$ is A001715.
\end{remark}


\subsection{Class B: $(123,\{0\},\{0, 1\})$}\label{c15}
The class is listed in Section \ref{sec:classB}.
\begin{proposition}
For any pattern $p$ in the class,
\[a_n(p)= n! - (n - 2)! (n - 2).\]
\end{proposition}
\begin{proof}
We consider only the pattern $p=(123,\{0\},\{0, 1\})$ here, the other cases being similar. Assume a permutation $\pi\in\perms_n$ has an occurrence of $p$. Then $\pi_1=1$, and for some $i$ and $j$, $i<j$, $\pi_i=2$ and $\pi_j>2$, and the remaining elements of $\pi$ are arbitrary. But since $\pi_1=1$ and $\pi_i=2$, $\pi_k>2$ for all $k>i$. Thus, the possible permutations are counted by $(n-2)(n-2)!$, since we can choose $i$ in $(n-2)$ ways.
\end{proof}
\begin{remark}
The sequence $c_n=(n-2)!(n-2)$ is A001563.
\end{remark}


\subsection{$(132,\{1, 2\},\{1, 2\})$ (1, 2, 5, 20, 102, 626, 4458,\ldots).} \label{c16}

These are the permutations without substrings of the form $(k,k+2,k+1)$. 

\begin{proposition}
Let $b_n=n!-a_n(132,\{1, 2\},\{1, 2\})$. Then
\[b_n=\sum_{k=1}^{\lfloor n/3\rfloor}(-1)^{k+1} (n - 2k)! \binom{n - 2k}{k}.\]
\end{proposition}
\begin{proof}
First note that each occurrence of the pattern, let us denote it by $p$, is of the form $(\pi_k,\pi_{k+1}+2,\pi_{k+2}+1)$ --- in the dot diagram it fills a 3 by 3 square. Let $c_{n,m,k}$ be the number of ways $k$ non-overlapping 3 by 3 squares can be put in an $n$ by $m$ rectangle. 

Then
\[c_{n,m, 1}=(n-2)(m-2),\] 
\[c_{n,m, 2}=2\sum_{i=3}^{n-3}\sum_{j=3}^{m-3}c_{i,j, 1},\]
and in general, 
\[c_{n,m,k}=k\sum_{i=3}^{n-3(k-1)}\sum_{j=3}^{m-3(k-1)}c_{i,j,k-1}.\] 
Now, 
\[c_{n,m, 2}=\frac{1}{2!}(n-4)(n-5)(m-4)(m-5),\]
\[c_{n,m, 3}=\frac{1}{3!}(n-6)(n-7)(n-8)(m-6)(m-7)(m-8),\]
and in general, 
\[c_{n,m,k}=\frac{1}{k!}\prod_{i=2k}^{3k-1}(n-i)(m-i).\] 
Further, 
\[(n-3k)!c_{n,n,k}=(n-2k)!\binom{n-2k}{k}.\]
Denote these numbers $d_{n,k}=(n-3k)!c_{n,n,k}$. 

Now, $d_{n,k}$ counts the number of permutations with at least $k$ occurrences of $p$, plus the number of permutations with at least $k+1$ occurrences, plus the number of permutations with at least $k+2$ occurrences, etc. Thus, 
\[b_n=d_{n, 1}-d_{n, 2}+d_{n, 3}-\cdots+(-1)^{\lfloor n/3\rfloor}d_{n,{\lfloor n/3\rfloor}}.\]
\end{proof}


\subsection{$(123,\{1, 2\},\{1, 2\})$ (A002628)}\label{c17}
Occurrence of the pattern means existence of a substring $k(k+1)(k+2)$. See the references in the OEIS for some classic results regarding problems of this kind.


\subsection{Class C: $(123,\{0\},\{0, 1, 2\})$}\label{c18}
See Section \ref{sec:classC} for a list of the members, modulo symmetries.
\begin{proposition}
For any pattern $p$ in the class, 
\[a_n(p)=n! - \frac{(n - 1)!}{2}.\]
\end{proposition}
\begin{proof}
Assume a permutation $\pi\in\perms_n$ has an occurrence of $(123,\{0\},\{0, 1, 2\})$. Then $\pi_1=1$, and for some $i$ and $j$, $i<j$, $\pi_i=2$ and $\pi_j=3$, and the remaining elements of $\pi$ are arbitrary. If $\pi_i=2$, there are $n-i$ possibilities for $j$. The total number of possibilities thus are
\[\left((n-2)+(n-3)+\cdots+1\right)(n-3)!=\frac{(n - 1)!}{2}.\]
\end{proof}
All patterns in the class can be enumerated in the same way.
\begin{remark}
The sequence $c_n=(n - 1)!/2$ is A001710.
\end{remark}


\subsection{Class D: $(123,\{0, 1\},\{0, 1\})$}\label{c19}
The class is listed in Section \ref{sec:classD}.
\begin{proposition}
For any pattern $p$ in the class, 
\[a_n(p)=n! - (n - 2)!.\]
\end{proposition}
\begin{proof}
Let $p=(123,\{0, 1\},\{0, 1\})$, the other cases being analogous. 
An occurrence of the pattern means that the permutation starts with $12$, and there are no further restrictions.
\end{proof}


\subsection{Class E: $(123,\{0, 1, 2\},\{0, 1, 2\})$}\label{c20}
This class is found in Section \ref{sec:classE}.
\begin{proposition}
For any pattern $p$ in the class, 
\[a_n(p)=n! - (n - 3)!.\]
\end{proposition}
\begin{proof}
All patterns in the class can be enumerated in the same way, so we only consider
 $p=(123,\{0, 1, 2\},\{0, 1, 2\})$. An occurrence of the pattern means that the permutation starts with $123$, with no restrictions on the remaining $n-3$ elements.
\end{proof}


\subsection{$(123,X,\{0, 1, 2, 3\})$ and $(132,X,\{0, 1, 2, 3\})$} \label{c21}
The requirements make it obvious that only one permutation, 123 or 132, depending on the case, contains the pattern, so we got the following result.
\begin{proposition}
\[a_n(123,X,\{0, 1, 2, 3\})=a_n(132,X,\{0, 1, 2, 3\})=n!-\mathbf{1}_{n=3}.\]
\end{proposition}

\appendix
\section{Some Wilf classes}
\subsection{Class A}\label{sec:classA}
{\small
\begin{align*}
&(123,\emptyset,\{0, 1\}), 
 (123,\emptyset,\{0, 3\}), 
 (123,\{0\},\{1, 2\}), 
 (123,\{0\},\{1, 3\}), 
 (123,\{0\},\{2, 3\}), 
 (123,\{1\},\{0, 2\}),\\ 
&(123,\{1\},\{0, 3\}), 
(123,\{1\},\{2, 3\}), 
 (132,\emptyset,\{0, 1\}), 
 (132,\emptyset,\{0, 3\}), 
 (132,\emptyset,\{2, 3\}), 
 (132,\{0\},\{1, 2\}),\\ 
& (132,\{0\},\{1, 3\}), 
 (132,\{0\},\{2, 3\}), 
 (132,\{1\},\{0, 1\}), 
 (132,\{1\},\{0, 2\}), 
 (132,\{1\},\{1, 3\}), 
 (132,\{1\},\{2, 3\}),\\ 
& (132,\{0, 1\},\{2\}), 
 (132,\{2\},\{0, 3\}), 
 (132,\{2\},\{1, 3\}),
 (132,\{0, 2\},\{3\}),
 (132,\{1, 2\},\{3\}), 
 (132,\{3\},\{0, 3\}),\\ 
& (132,\{3\},\{1, 3\}).
\end{align*}
}
\subsection{Class B}\label{sec:classB}
{\small \begin{align*}
&(123,	\{0\},	\{0, 1\}),
(123,	\{0\},	\{0, 3\}),
(123,	\{1\},	\{0, 1\}),
(123,	\{0, 1\},	\{0, 2\}),
(123,	\{0, 1\},	\{1, 2\}),\\
&(123,	\{0, 1\},	\{0, 3\}),
(123,	\{0, 1\},	\{1, 3\}),
(123,	\{0, 1\},	\{2, 3\}),
(123,	\{0, 2\},	\{1, 2\}),
(123,	\{0, 2\},	\{0, 3\}),\\
&(123,	\{0, 2\},	\{1, 3\}),
(123,	\{1, 2\},	\{0, 3\}),
(132,	\{0\},	\{0, 1\}),
(132,	\{0\},	\{0, 3\}),
(132,	\{1\},	\{0, 3\}),\\
&(132,	\{0, 1\},	\{0, 1\}),
(132,	\{0, 1\},	\{0, 2\}),
(132,	\{0, 1\},	\{1, 2\}),
(132,	\{0, 1\},	\{3\}),
(132,	\{0, 1\},	\{1, 3\}),\\
&(132,	\{0, 1\},	\{2, 3\}),
(132,	\{2\},	\{2, 3\}),
(132,	\{0, 2\},	\{1, 2\}),
(132,	\{0, 2\},	\{0, 3\}),
(132,	\{0, 2\},	\{1, 3\}),\\
&(132,	\{0, 2\},	\{2, 3\}),
(132,	\{1, 2\},	\{0, 3\}),
(132,	\{1, 2\},	\{1, 3\}),
(132,	\{1, 2\},	\{2, 3\}),
(132,	\{3\},	\{2, 3\}),\\
&(132,	\{0, 3\},	\{0, 3\}),
(132,	\{0, 3\},	\{1, 3\}),
(132,	\{0, 3\},	\{2, 3\}),
(132,	\{1, 3\},	\{1, 3\}),
(132,	\{1, 3\},	\{2, 3\}).
\end{align*}}
\subsection{Class C}\label{sec:classC}
{\small\begin{align*}
&(123,	\{0\},	\{0, 1, 2\}),
(123,	\{0\},	\{0, 1, 3\}),
(123,	\{0\},	\{0, 2, 3\}),
(123,	\{0\},	\{1, 2, 3\}),
(123,	\{1\},	\{0, 1, 2\}),\\
&(123,	\{1\},	\{0, 1, 3\}),
(123,	\{1\},	\{0, 2, 3\}),
(123,	\{1\},	\{1, 2, 3\}),
(132,	\{0\},	\{0, 1, 2\}),
(132,	\{0\},	\{0, 1, 3\}),\\
&(132,	\{0\},	\{0, 2, 3\}),
(132,	\{0\},	\{1, 2, 3\}),
(132,	\{1\},	\{0, 1, 2\}),
(132,	\{1\},	\{0, 1, 3\}),
(132,	\{1\},	\{0, 2, 3\}),\\
&(132,	\{1\},	\{1, 2, 3\}),
(132,	\{2\},	\{0, 1, 2\}),
(132,	\{2\},	\{0, 1, 3\}),
(132,	\{2\},	\{0, 2, 3\}),
(132,	\{2\},	\{1, 2, 3\}),\\
&(132,	\{0, 1, 2\},	\{3\}),
(132,	\{3\},	\{0, 1, 3\}),
(132,	\{3\},	\{0, 2, 3\}),
(132,	\{3\},	\{1, 2, 3\}).
\end{align*}}
\subsection{Class D}\label{sec:classD}
{\small \begin{align*}
&(123,\{0, 1\},\{0, 1\}),
(123,\{0, 1\},\{0, 1, 2\}),
(123,\{0, 1\},\{0, 1, 3\}),
(123,\{0, 1\},\{0, 2, 3\}),\\
&(123,\{0, 1\},\{1, 2, 3\}),
(123,\{0, 2\},\{0, 1, 2\}),
(123,\{0, 2\},\{0, 1, 3\}),
(123,\{0, 2\},\{0, 2, 3\}),\\
&(123,\{0, 2\},\{1, 2, 3\}),
(123,\{1, 2\},\{0, 1, 2\}),
(123,\{1, 2\},\{0, 1, 3\}),
(123,\{0, 1, 2\},\{0, 3\}),\\
&(123,\{0, 3\},\{0, 3\}),
(123,\{0, 3\},\{0, 1, 3\}),
(132,\{0, 1\},\{0, 1, 2\}),
(132,\{0, 1\},\{0, 3\}),\\
&(132,\{0, 1\},\{0, 1, 3\}),
(132,\{0, 1\},\{0, 2, 3\}),
(132,\{0, 1\},\{1, 2, 3\}),
(132,\{0, 2\},\{0, 1, 2\}),\\
&(132,\{0, 2\},\{0, 1, 3\}),
(132,\{0, 2\},\{0, 2, 3\}),
(132,\{0, 2\},\{1, 2, 3\}),
(132,\{1, 2\},\{0, 1, 2\}),\\
&(132,\{1, 2\},\{0, 1, 3\}),
(132,\{1, 2\},\{0, 2, 3\}),
(132,\{1, 2\},\{1, 2, 3\}),
(132,\{0, 1, 2\},\{0, 3\}),\\
&(132,\{0, 1, 2\},\{1, 3\}),
(132,\{0, 1, 2\},\{2, 3\}),
(132,\{0, 3\},\{0, 1, 3\}),
(132,\{0, 3\},\{0, 2, 3\}),\\
&(132,\{0, 3\},\{1, 2, 3\}),
(132,\{1, 3\},\{0, 1, 3\}),
(132,\{1, 3\},\{0, 2, 3\}),
(132,\{1, 3\},\{1, 2, 3\}),\\
&(132,\{0, 1, 3\},\{2, 3\}),
(132,\{2, 3\},\{2, 3\}),
(132,\{2, 3\},\{0, 2, 3\}),
(132,\{2, 3\},\{1, 2, 3\}).
\end{align*}}
\subsection{Class E}\label{sec:classE}
{\small \begin{align*}
&(123,\{0, 1, 2\},\{0, 1, 2\}),
(123,\{0, 1, 2\},\{0, 1, 3\}),
(123,\{0, 1, 2\},\{0, 2, 3\}),
(123,\{0, 1, 2\},\{1, 2, 3\}),\\
&(123,\{0, 1, 3\},\{0, 1, 3\}),
(123,\{0, 1, 3\},\{0, 2, 3\}),
(132,\{0, 1, 2\},\{0, 1, 2\}),
(132,\{0, 1, 2\},\{0, 1, 3\}),\\
&(132,\{0, 1, 2\},\{0, 2, 3\}),
(132,\{0, 1, 2\},\{1, 2, 3\}),
(132,\{0, 1, 3\},\{0, 1, 3\}),
(132,\{0, 1, 3\},\{0, 2, 3\}),\\
&(132,\{0, 1, 3\},\{1, 2, 3\}),
(132,\{0, 2, 3\},\{0, 2, 3\}),
(132,\{0, 2, 3\},\{1, 2, 3\}),
(132,\{1, 2, 3\},\{1, 2, 3\}).
\end{align*}}

\end{document}